\documentclass[smallextended]{article}       

\usepackage{authblk}
\usepackage{lipsum}
\usepackage{amsmath}
\usepackage{amssymb}
\usepackage{todonotes}
\usepackage{graphicx}
\usepackage{pythonhighlight}
\usepackage{multicol}
\usepackage{placeins} 

\makeatletter
\let\cl@chapter\undefined
\makeatletter
\usepackage{cleveref}
\usepackage{mathtools}

\newcommand{\RR}{\mathbb{R}}

\newcommand{\Vh}{\mathcal{V}^h}
\newcommand{\VDh}{\mathcal{V}^h_{\mathrm{D}}}
\newcommand{\HD}{\mathcal{H}^2_{\mathrm{D}}(\Omega)}
\newcommand{\HH}{\mathcal{H}^2 (\Omega)}

\newcommand{\OD}{\Omega_{\mathrm{D}}}
\newcommand{\Oe}{\Omega_{\mathrm{e}}}

\newcommand{\Om}{\Omega_{\mathrm{m}}}

\newcommand{\Ge}{\Gamma_{\mathrm{e}}}

\newcommand{\uh}{u^h}

\newcommand{\ndof}{n_{\mathrm{dof}}}

\newcommand{\dx}{\;\mathrm{d}x}

\newcommand{\trans}{^{\mathrm{T}}}

\crefname{equation}{}{}
\Crefname{equation}{Equation}{Equations}
\crefname{section}{Sec.}{Secs.}
\Crefname{section}{Sec.}{Section}
\crefname{figure}{Fig.}{figures}
\Crefname{figure}{Fig.}{Figures}
\usepackage[utf8]{inputenc}
\usepackage[T1]{fontenc}

\usepackage[english]{babel}

\title{Modeling a Smooth Surface by a Constrained Biharmonic Equation with Application in Soil
  Science\thanks{Supported by NORPART-2021/10167 Ethiopian Norwegian Network in Computational
    Mathematics (ENNCoMat).}
}
\author[1]{Samson~Seifu~Bekele} 
\author[1]{Maregnesh~Mechal~Wolde}
\author[2]{Claus~Führer}
\author[3]{Nils-Otto~Kitterød}
\author[5]{Anne~Kværnø}

\affil[1]{Department of Mathematics, Hawassa University, Hawassa, Ethiopia}
\affil[2]{Center for Mathematical Sciences,
           Lund University (LU), 
           Lund, Sweden}
\affil[3]{
           Soil and Water Science Section, 
           Norwegian University of Life Sciences (NMBU) ,
           {\AA}s, Norway}  
\affil[4]{
           Department of Mathematical Sciences,
           Norwegian University of Science and Technology (NTNU), Trondheim, Norway}

\date{}

\begin{document}
\maketitle

\begin{abstract}
This paper presents a method for mathematical modelling of surfaces conditioned on empirical data. It is based on solving
a discrete biharmonic equation over a domain with given inner point and inner curve data. The inner 
curve data is used to model boundary values while the inner point data is used for 
modeling a load vector with the goal 
to generate a smooth surface. The construction of boundary data is an ill-posed problem, for which a 
special regularization approach is suggested.

The method is designed for surface construction problems with a very limited amount of measured data. In the paper 
we apply the method by using empirical data of soil thickness and geological maps indicating  
exposed bedrock regions.


\end{abstract}

\section{Introduction} \label{sec:intro} 
This paper 
deals with constructing a surface over a  2D domain from boundary data, 
point data and data on curves inside the domain. Surface reconstruction 
is an important topic in many applications and related tools are often bound to 
specific mathematical models for the application at hand. Smoothness 
requirements, free parameters for adapting the final shape to given 
data or the need for compensation for lack of data are often criteria 
for the design of a mathematical surface model. 

There are various ways 
to choose surface models. In computer aided graphical design (CAGD) the 
models are based on function classes like splines and NURBS
\cite{NURBS2000}. In CAGD given data points usually determine the free function coefficients. If the number of 
given data points fits the dimension of the function space, the surface 
construction is a 2D-interpolation problem. In the case that more 
data is given, surface construction leads to a fitting problem. Spline 
or NURBS function spaces are chosen such that data change affects 
the surface only locally. 

In the current study we consider situations where much less 
data is available. This is the case in the guiding example from soil science
where soil thickness has to be shaped based on very little 
data from occasional drilling down to the bedrock. The approach we consider for this case is based on filling the 
missing data by using a smooth surface model generated as a solution of 
a discretized
partial differential equation (PDE). 
Typically, load functions  and boundary 
values are used to shape the surface. In this paper we propose a method 
for the case when additionally point data or data from interior 
curves is available and when 
boundary data is not completely known. A PDE-based 
model allows the use of 
well-established numerical tools and software and is suited for irregular 
domains with non uniformly distributed measurements. An example of such an approach, 
based on the Poisson equation, is found in \cite{kazhdan2006}. 
A nonlinear PDE-based method, as discussed in \cite{claisse2011nonlinear}, utilizes the level set approach, where a higher-dimensional scalar function, known as the level set function, evolves according to the PDE. This function implicitly represents the surface through its zero level curve, ensuring 
the reconstruction of a smooth and regular surface that closely aligns with the sampled data points.

Inspired by plate bending problems in solid mechanics, we present in 
this paper the use of 
the discretized biharmonic equation to model surfaces and apply the 
virtual element 
method to define the discrete model space. The approach is demonstrated 
through the example of determining soil thickness over a specific 
geographical region. The biharmonic equation requires boundary data 
like position, derivative information, and load data over the entire 
domain. In surface reconstruction, this data is often not directly 
available and is instead used as modeling parameters to create surfaces 
with a desired shape, such as interpolating specific points, containing 
given lines, or meeting certain constraints. The paper takes up the 
following aspects: the formulation of the mathematical problem, the presentation of a guiding
application, modeling alternatives, solution methods and a brief summary of the results.

\subsection{Problem description} \label{sec:problem}
Our main objective is to construct the surface under consideration as the solution of a certain discretized biharmonic equation. 
The domain, the load vector and the boundary conditions of this biharmonic equation serve as modeling parameters to shape this surface 
in a way that given measurement data is taken into account. For this end we consider the following problem  on a closed domain $\Omega \subset  \RR^2$:
\begin{subequations}\label{eq:ourproblem}
\begin{equation} 
\Delta^{2}u = q  \;\; \mathrm{in} \;\; \Omega \label{eq:ourproblem_a} 
\end{equation}
with boundary conditions
\begin{align}
  &&u                       & = f \;\; & \text{and} &
  & \partial _n u           & = h \;\; & \text{on} \;\; \partial \Omega_{\mathrm{D}}, &
  & \label{eq:ourproblem_b}&& \\
  &&\Delta u                & = 0 \;\; & \text{and} &
  & \partial _n (\Delta u)  & = 0 \;\; & \text{on} \;\; \partial \Omega_{\mathrm{N}}. &
  & \label{eq:ourproblem_c} &&
\end{align}
\end{subequations}
The boundary
$\partial \Omega = \partial \Omega_{\mathrm{D}} \cup \partial \Omega_{\mathrm{N}}
$ is divided into a part $\partial \Omega_{\mathrm{D}} $ for Dirichlet boundary conditions and
another $\partial \Omega_{\mathrm{N}} $ for Neumann boundary 
conditions. The right hand side function $q$
will be referred to as the load function.

Given the Sobolev space $\HH$ and the subspace
\[
\HD= \{v \in \HH \colon  v|_{\partial \OD} = 0, (\partial _n v)|_{\partial \OD}=0 \},
\]
the variational formulation of \cref{eq:ourproblem} is: \\[2ex]
Find $u \in \HH$ satisfying the boundary conditions  \cref{eq:ourproblem_b,eq:ourproblem_c} and
\begin{equation} \label{eq:variationalproblem}
\int _{\Omega} \Delta u \Delta v \dx=\int _{\Omega}qv\dx, \;\; \forall v \in \HD .
\end{equation}

The challenge is then to decide the load function $q$ and the boundary functions $f$ and $h$ based on
point data given either as equality or inequality constraints:
\begin{subequations}
  \label{eq:ourproblem2}
\begin{align}
  u(\xi_i) & = d_i \;\; \mathrm{for} \;\; i \in \mathcal{E}\label{eq:ourproblem2e}, \\
  u(\xi_i) & \geq d_i\;\; \mathrm{for} \;\;  i \in \mathcal{I} \label{eq:ourproblem2i},
\end{align}
\end{subequations}
where $\xi_i \in \Omega$ and  $\mathcal{E}$ and $\mathcal{I}$
are finite sets of indices. 
Furthermore, the surface values may be known at a set of curves in the domain:
\begin{align}
  u &=  u_{\mathrm{e}} \;\; \mathrm{on} \;\;   \Gamma_{\mathrm{e}} \subset \Omega.  \label{eq:ourproblem1}
\end{align} 

In \cref{sec:loadvector} we describe a procedure for using point data \cref{eq:ourproblem2} to obtain parameters
of a load function $q$ model. In the case, when the resulting surface is of importance only in a small subarea of the domain in which all point data is available, the effect of
the boundary values can be neglected.   
For the other case,  we describe in \cref{sec:boundary} how missing Dirichlet boundary
data \cref{eq:ourproblem_b} can be constructed from information at a set of curves, as expressed in 
\cref{eq:ourproblem1}. Both procedures are applied to an example from soil science, described in the
next section.


\subsection{Organization of the paper}
The topic of this paper is to cover two aspects of surface construction by an PDE approach: 
In \cref{sec:loadvector} we focus on the load vector as a modeling parameter. Inner point constraints are used 
to optimally fit these parameters while in \cref{sec:boundary} given curve data is used to determine 
values for Dirichlet boundary conditions of the desired surface.
The domains can be geometrically complex in general, while the approach of including boundary data as modeling parameters 
allows for simplifications of the domain to simple rectangles.
Both aspects of surface construction are inherent in the surface construction problem from geoscience introduced
in \cref{sec:bedrock}  and used to demonstrated the modeling approaches in \cref{sec:loadvector,sec:boundary}. 
The corresponding software tools are presented in \cref{sec:software}.

\section{Case study: Soil thickness and bedrock surface over a region in Norway} \label{sec:bedrock}
\begin{figure}[h]
  \centering
  \includegraphics[width=0.5\textwidth]{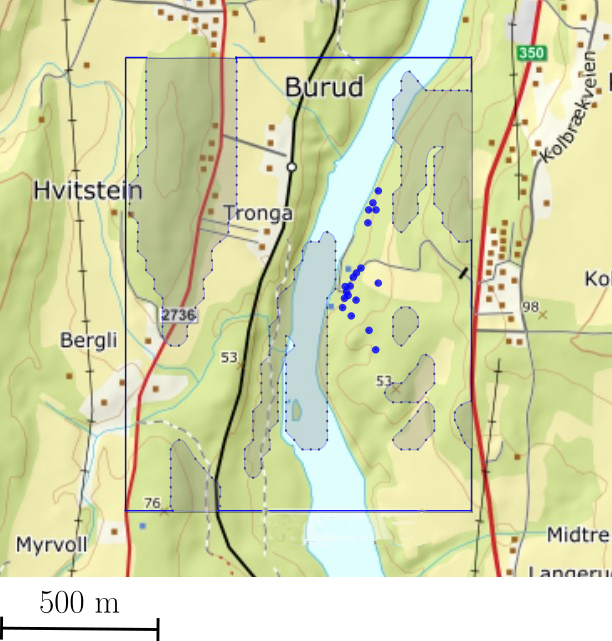}
  \caption{A part of the Norwegian region Øvre Eiker with the subregion 
  marked for which the soil thickness is determined, \cite{norgeskarta} 
  region:  EU89 UTM33 (6639960\,N, 212750\,E). The position of the wells (point data) is indicated by blue circles  and the exposed bedrock $\Oe$ is indicated by slightly grayed areas.}
\label{fig:areas}
\end{figure}

Consider a given rectangular geographical area  $\Omega_{\mathrm{m}}$ with a subarea, $\Omega_{\textrm c}$,  covered by
soil. In the remaining part of the area, $\Omega_{\mathrm{e}} := \Omega_{\mathrm{m}} \backslash \Omega_{\mathrm{c}}$, the bedrock is exposed, cf.~\cref{fig:areas}.
Correspondingly, let $u_{\mathrm{m}}$ be the terrain surface elevation,  $u_{\mathrm{s}}$ the soil thickness and
$u_{\mathrm{b}}$ the bedrock surface elevation, thus
\begin{equation}
\label{eq:sumofboth}
   u_{\mathrm{b}} = u_{\mathrm{m}} - u_{\mathrm{s}}, \quad \text{in} \quad \Omega_{\mathrm{m}}, 
\end{equation}
with $u_{\mathrm{s}}=0$ in  the exposed area $ \Omega_{\mathrm{e}}$.

We access the terrain surface data from the database \cite{kartverket2024} and the
location of the exposed bedrock from  \cite{ngu2024a}.
Information about soil thickness and bedrock morphology is relevant for many different applications (\emph{viz.} agriculture; 
hydrology, structural engineering). In this work, soil thickness will be modelled based on local point information. Bedrock morphology will be modelled based on information about the elevation of the exposed bedrock area, in particular on the intersections between the exposed and the covered regions.

The point data that will be used in  \cref{eq:ourproblem2} are the locations and depths of 
wells drilled in $\Omega_{\mathrm{c}}$. It can be
retrieved from a public database \cite{ngu2024b}.
When a well is drilled down to the bedrock, the soil thickness $u_s$  is known at that point,
satisfying the equality constraint \cref{eq:ourproblem2e}. The inequality constraints 
\cref{eq:ourproblem2i} refer in this context to wells that do not reach 
entirely down to 
the bedrock. These data will be used in \cref{sec:loadvector} to construct a surface representing the soil thickness.

\begin{figure}[h]
\begin{center}
\includegraphics[width=0.75\textwidth]{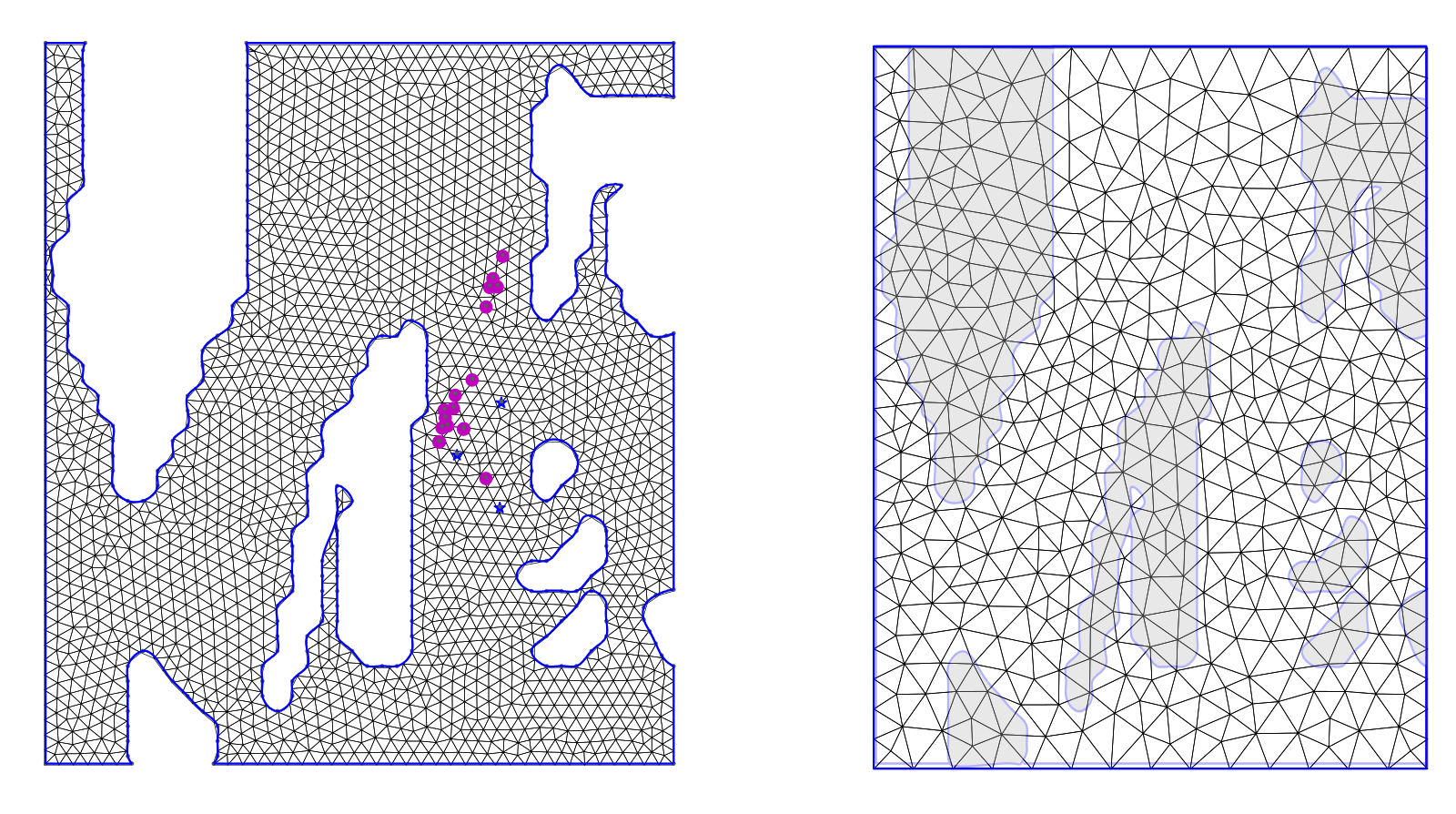} \hfill
\caption{Two alternatives to define the underlying domain together with an automatically generated grid.
{\it Left:} An irregular domain defined by a rectangular map segment where regions with exposed bedrock are excluded. 
The  location of wells is marked by magenta circles and blue stars. {\it Right:} The same rectangular map segment including regions with  
exposed bedrock (shaded).}
\label{fig:meshes}
\end{center}
\end{figure}

The terrain surface $u_{\mathrm{m}}$ is known over the entire domain $\Omega_{\mathrm{m}}$, so are the
Dirichlet boundary conditions $f$ in \cref{eq:ourproblem_b} on the boundaries between the exposed and 
the covered regions in \cref{sec:loadvector}.
In \cref{sec:boundary},
these boundaries will act as the set of curves $\Gamma_{\mathrm{e}}$ in \cref{eq:ourproblem1}, and will be used to retrieve
the unknown Dirichlet boundary values at the part of the boundaries 
covered by soil. Hence we consider two different domains to formulate the problem: 
In \cref{sec:boundary} the domain $\Omega$  is completely defined by 
the rectangular map section $\Omega_{\mathrm{m}}$. Here, 
the boundary is given by the edges of the rectangle. In \cref{sec:loadvector}
the domain is defined as $\Omega:=\Omega_{\mathrm{c}}$. It is derived from the map section by cutting out regions 
with exposed bedrock. The edges of those parts together with the 
edges of the rectangular map form the boundary of  
$\Omega_{\mathrm{c}}$. In \cref{fig:meshes} these two domains are shown together. One sees how the more 
complex domain geometry when the exposed areas are removed from the domain, affects the mesh.

The aim is to reproduce a smooth large-scale trend of the soil thickness $u_s$ and the bedrock
surface $u_b$.
Different approaches of modelling the soil thickness and the bedrock 
surface have been suggested before. An example is 
modeled  hill slope sediment thickness using analytical
solution of differential equations, viewing sediment transport as a 
diffusive process, \cite{RemDiet2014}.  For a Poisson equation 
based surface model to describe soil thickness we refer to  
\cite{kitterod2017estimating}.
Further references to this problem can be found in
\cite{kitterod2021estimation}.


\section{Modelling the load function from point data}  \label{sec:loadvector}

In this section, we describe a procedure for using point data \cref{eq:ourproblem2} to model a load
function $q$ in the biharmonic equation \cref{eq:ourproblem}. 
In this case, the boundary conditions in
\cref{eq:ourproblem_b} are assumed to be known, such that the smooth surface is given as the
solution $u$ of the discrete biharmonic equation. The procedure is then applied to the soil thickness
problem described in \cref{sec:bedrock}. 

The general idea can be described as follows: Let $q$ be a weighted sum of load functions $q_j$ representing
each data point, thus  $q = \sum_{j=1}^{n_q} p_j q_j$. In the current work, these are chosen in relation to the locations $\xi_j$ of the point data as
\begin{equation}\label{eq:Gaussian}
q_j(x) = \frac{1}{\sigma_j \sqrt{2 \pi}}e^{- \frac{1}{2}(\frac{x- \xi_j}{\sigma_j})^2}, \qquad j\in
\mathcal{E}\cup \mathcal{I}.
\end{equation}
These functions respect local information,
and the choice of the localization parameters $\sigma_j$ controls the region of their 
influence. Let $u_j$ be the solutions of the
biharmonic equations $\Delta^2 u_j = q_j$ with homogeneous Dirichlet boundary conditions for each of
the basis functions $q_j$, $j=1,\dotsc,n_q$. 
For a given function $q = \sum_{j=1}^{n_q} p_j q_j$ with weights $p=(p_1,\dotsc,p_{n_q})$ the solution of the biharmonic equation
\cref{eq:ourproblem} is   
\begin{equation} \label{eq:uall}
  u(x;p) = \sum_{j=1}^{n_q} p_j u_j(x) + u_{\mathrm{bc}}(x),
\end{equation}
where $u_\mathrm{bc}$ is the solution of $\Delta^2 u_\mathrm{bc}=0$ respecting the given boundary
conditions. The challenge is then to decide on the weight vector $p$. In this work, $p$ is
found by solving the following minimization problem:

\begin{subequations} 
  \label{eq:minimization}
\begin{equation} \label{eq:minimization0}
  \text{min}_{p \in \mathcal{B}} \quad  \sum_{i = 1}^{n_q}(d_i - u({\xi_i;p}))^2
\end{equation}
subject to 
\begin{align}
  u(\xi_i;p) & = d_i \;\; \mathrm{for} \;\;  i \in \mathcal{E}, \label{eq:minimizatione} \\
  u(\xi_i;p) & \geq d_i \;\; \mathrm{for} \;\;  i \in \mathcal{I}, \label{eq:minimizationi}
\end{align}
\end{subequations}
where $\mathcal{B} = [a_1,b_1] \times \cdots \times [a_{n_q}, b_{n_q}]$, with $[a_i,b_i]$ bounding
the weight $p_i$, respectively. A suitable choice of bounds has to be considered to ensure that the surface generated by the biharmonic
equation \cref{eq:ourproblem} interpolates the specified equality  data points \cref{eq:minimizatione},
without necessarily interpolating the inequality points, and at the same time to control
the smoothness of the
surface to suppress possible unwanted oscillations.  Without bounds, the solution will interpolate all data points, which is not the
intention. The exact choice of the bounds depends on the actual application. 

\subsection{Case study: Modelling soil thickness}\label{Sec:experiment}
\begin{figure}[h]
\begin{center} 

\includegraphics[width=0.4\textwidth]{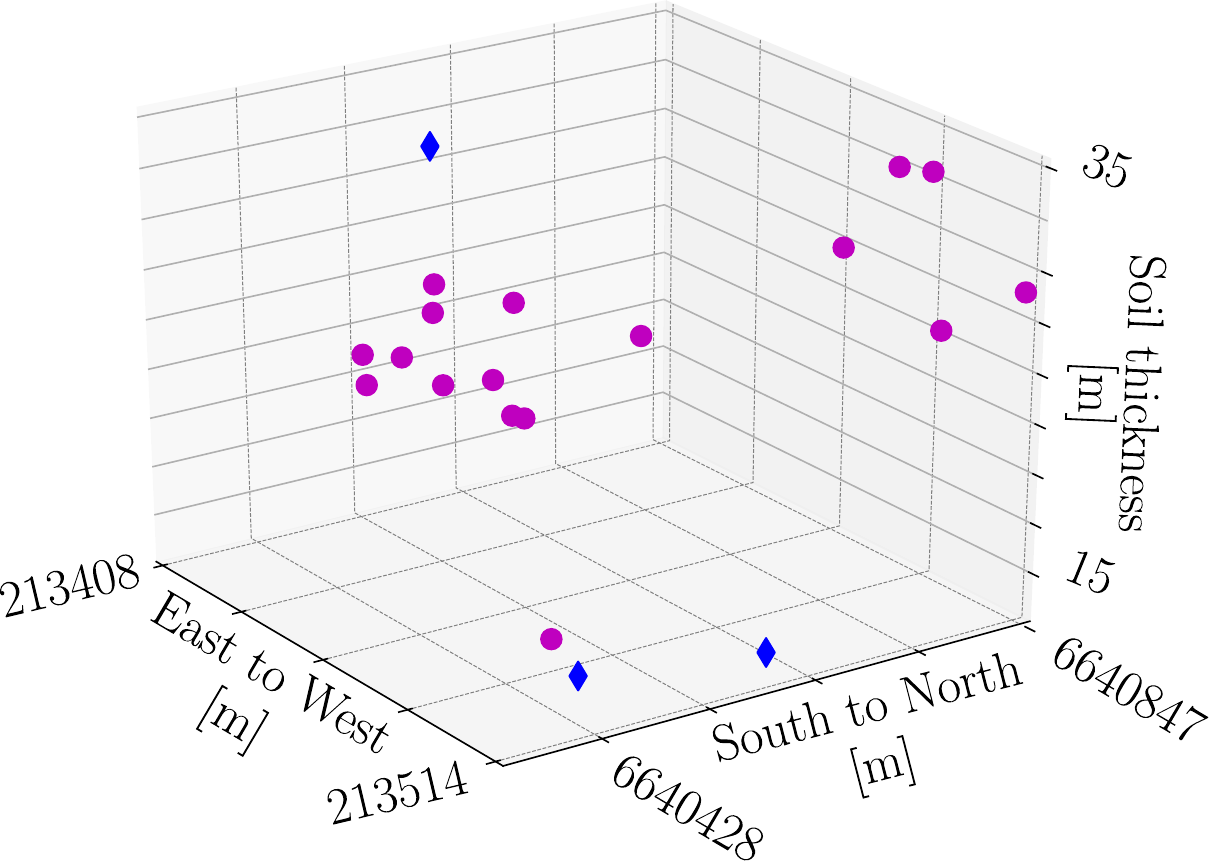}

\caption{Distribution of the point data. Exact soil thickness (equality constraints) $\color{blue}
  \blacklozenge$. Lower bound of soil thickness (inequality constraints) 
  $\color{magenta} \bullet$. }
\label{fig:M1}
\end{center}

\end{figure}

We will now apply the procedure described above to derive a surface mimicking the soil thickness $u_s$
by using information from the wells as point data, as described in \cref{sec:bedrock}. 
The computational domain is restricted to the area covered by soil, thus 
$\Omega=\Omega_{\mathrm{c}}$. The solution $u$ represents the soil thickness $u_s$. At the
boundaries between the exposed and the covered regions, the soil thickness is obviously
zero. Thus, these boundaries form the Dirichlet boundary $\partial \Omega_D$, and $u$ and $\partial_n u$ are set to be
zero.  On the remaining parts of the boundary
information is missing; these boundaries are chosen as zero Neumann boundaries $\partial \Omega_N$ \cref{eq:ourproblem_c}. 
The data set consists of information from twenty wells. The distribution of the wells is visualized in \cref{fig:M1}. Only three of the wells are drilled down to the
bedrock, here the soil thickness $d_i$, $i\in\mathcal{E}$ is known. This is expressed by equality
constraints \cref{eq:ourproblem2e},  the blue diamonds in \cref{fig:M1}.  For the remaining wells, the soil thickness is larger than the
well depth $d_i$, $i\in\mathcal{I}$, so information from them is expressed by inequality
constraints \cref{eq:ourproblem2i}. The corresponding well depths are represented by magenta dots in
\cref{fig:M1}.

\begin{figure}[h]
\begin{center}
  \includegraphics[width=0.3\textwidth]{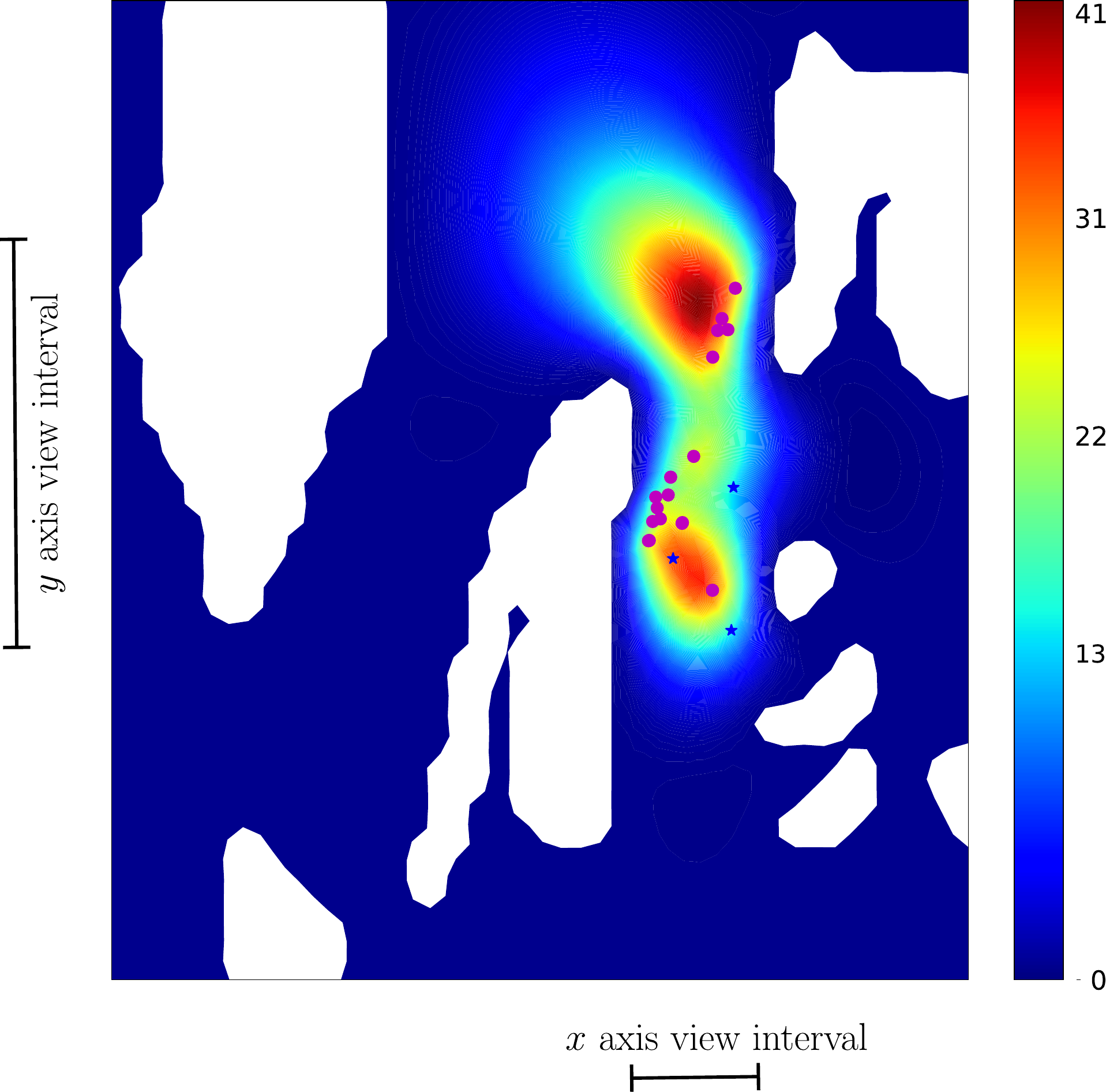} \hfill
  \includegraphics[width=0.32\textwidth]{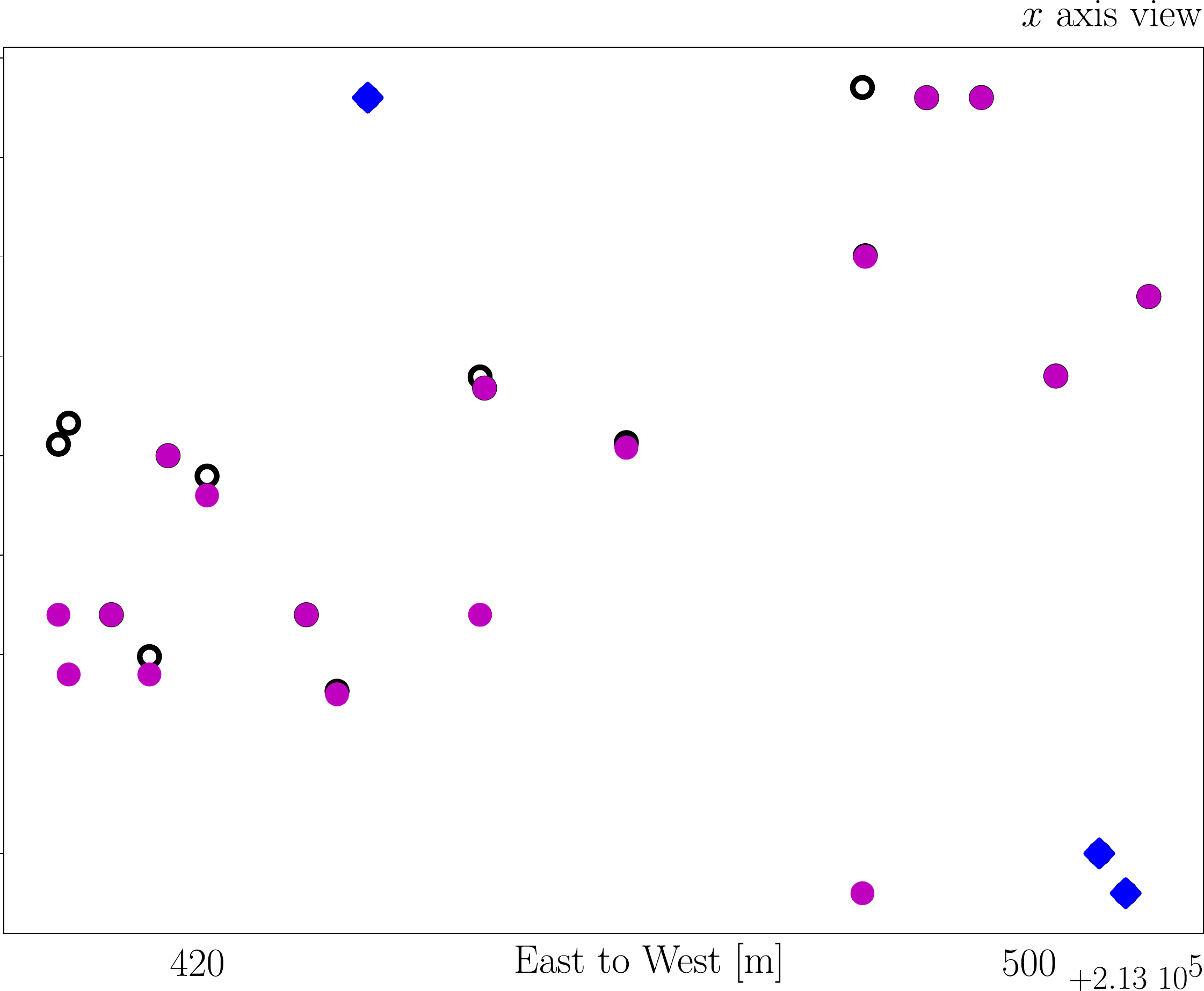} 
  \includegraphics[width=0.34\textwidth]{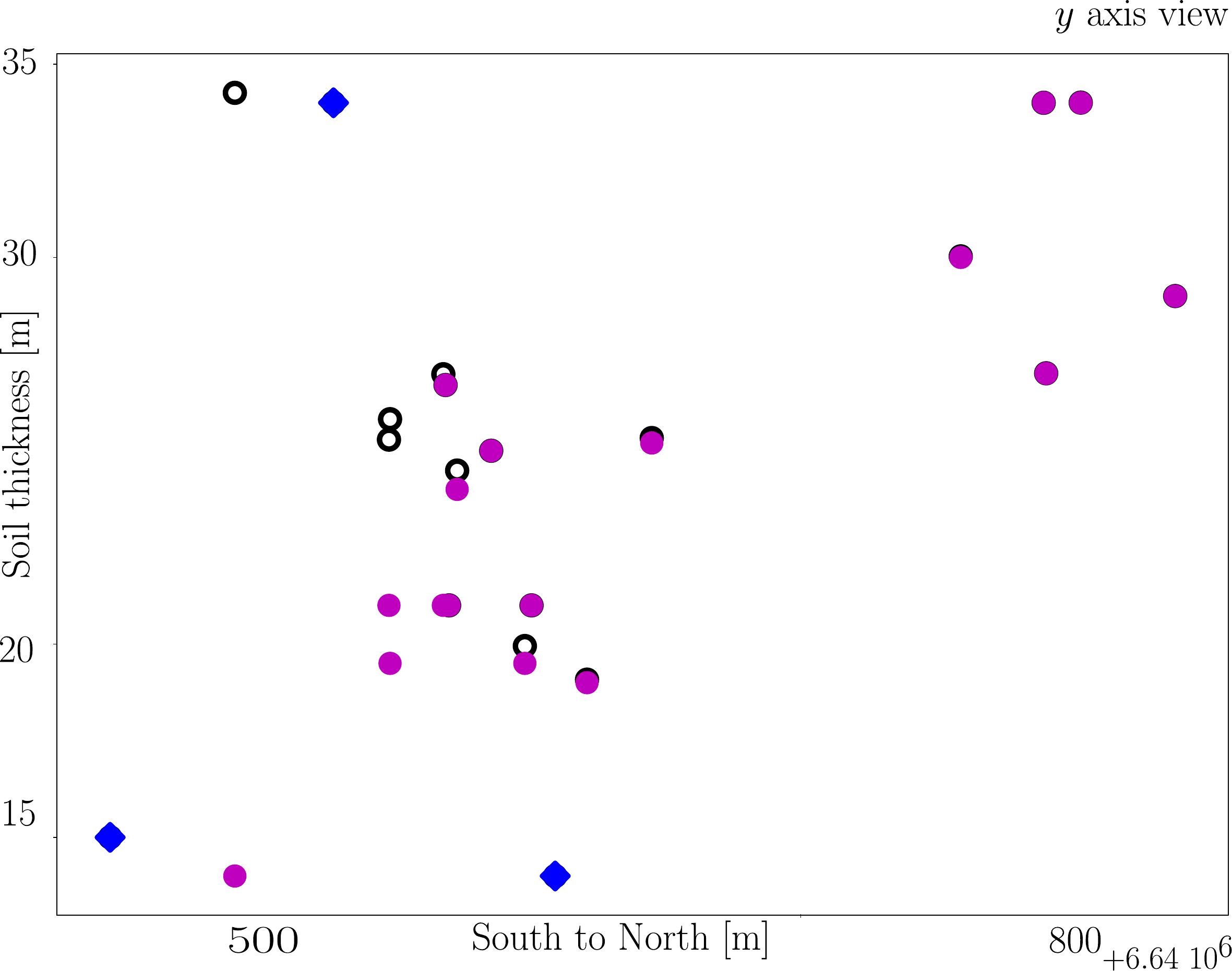} 
\caption{{\it Left:} Surface produced from the well data. {\it Right:} Measured and simulated point data, seen from the north-south axis and seen from the east-west axis: Given soil thickness $\color{blue}
  \blacklozenge$, well depths (lower bound of soil thickness)  $\color{magenta} \bullet$, simulated
  soil thickness {\large $\circ$}. }
\label{fig:M2}
\end{center}
\end{figure}
\FloatBarrier


One example of a surface constructed from the well data is depicted in \cref{fig:M2}, to the left. In this example, the localization
parameters in \cref{eq:Gaussian} are set to  $\sigma_i=1.5$  for all points $\xi_i$.  Suitable
bounds for 
the weights are chosen experimentally to 
$p_i \in [0,2]$ for $i\in\mathcal{E}$ and $p_i \in [-1,1]$ for
$i\in \mathcal{I}$.  The discrete  biharmonic equation is
constructed by the Virtual Element Method (VEM) on the grid given in \cref{fig:meshes}, to the
left. 
All the measured values $(\xi_i,d_i)$ as well as the simulated values
$(\xi_i,u(\xi_i,p))$ are presented in \cref{fig:M2}, to the right,
confirming that all the constraints are satisfied.

This approach uses the local information from the wells to give a rough estimate on
the soil thickness in a small subdomain that contains all wells. 
The smaller this subdomain the less
is the influence of the outer boundary data, that normally is unknown. 
In contrast, the next section suggests a technique to estimate boundary data from inner curves. These are the
boundary lines between the exposed and the unexposed bedrock, which will be used to estimate the bedrock topography over a larger area.


\section{Modelling outer boundaries from inner curve data} \label{sec:boundary}
We consider now the case of constructing boundary conditions from given 
data on a set of curves $\Ge \subset \Omega$
while data on the boundary is not completely available.  In particular 
we consider \cref{eq:ourproblem} as a pure Dirichlet problem with $f$ and $h$ being unknown 
on the entire boundary $\partial \Omega_\mathrm{D}$ or parts  of it. 
Again, the application presented in \cref{sec:bedrock} serves as a case study after a solution method has been presented.

The problem is formulated as follows:
\begin{subequations}   \label{eq:constrweekform}
Given two functions $u_{\mathrm{e}_1}$ on 
$\Gamma_{\mathrm{e}_1}  \subseteq \Ge$ and $u_{\mathrm{e}_2}$ on 
$\Gamma_{\mathrm{e}_2} \subseteq \Ge$ 
find $u \in \HH$  such that
\begin{eqnarray} \label{eq:variationalproblem_2}
\int _{\Omega} \Delta u \Delta v \dx   & =  & \int_{\Omega} q v \dx, \: 
\;\forall v \in \HD \\ \label{eq: var_pro}
u &= & u_{\mathrm{e}_1}  \; \text{on} \; \Gamma_{\mathrm{e}_1} \\ \label{eq:Interior_cond_1}
\partial _n u &= & u_{\mathrm{e}_2} \; \text{on} \; \Gamma_{\mathrm{e}_2}.\label{eq:Interior_cond_2}
\end{eqnarray}
\end{subequations}

This formulation  is too general to allow a 
statement about its well-posedness. Indeed, inconsistent choices of the sets 
$\Gamma_{\mathrm{e}_1}$ and $\Gamma_{\mathrm{e}_2}$ might result in a 
problem without a solution. But it may  still
be possible to find a solution of its
discretized counterpart, which will be used in 
\cref{sec:discrete_inverse}.

On the other hand, given that a  solution of \cref{eq:constrweekform} 
exists, it 
defines the unknown boundary functions by  
\[
f:=\left . u \right |_{\partial \Omega _D} \text{ and } h:= \left . \partial _n u \right |_{\partial \Omega _D}.
\]

Note that, the well posed problem \cref{eq:ourproblem} is a special case of \cref{eq:constrweekform} with $\Gamma_{\mathrm{e}_1} = \Gamma_{\mathrm{e}_2} = \partial \Omega _D$.

\subsection{The discrete inverse problem}
\label{sec:discrete_inverse}

A discrete formulation of problem \cref{eq:constrweekform}
can be obtained by a Ritz ansatz with an $\ndof$-dimensional function space 
$\Vh=\text{span}\{\phi_1,\cdots,\phi_{\ndof}\}\subset \HH$ and a subspace
$\VDh = \{ \phi \in \Vh : \phi = 0 \text{ and } \partial_n \phi = 0 
\text{  on  } \partial \OD \} \subset \HD$.
Typically, the basis functions are functions with local support as used by the finite element
method or the virtual element method. 
We denote the set of the indexes of the basis functions $\phi_j \in 
\VDh$ by $\mathcal{J}$ and their total number by $| \mathcal{J}|$.
Furthermore we choose  
a discrete set of points on the inner curves, $\{\xi_l\}_{l=1}^n  
\subset \Ge$.
The discrete solution of \cref{eq:constrweekform} is 
\begin{equation*}
  \uh := \sum_{j=1}^{\ndof} c_j \phi_j 
\end{equation*}
with coefficients $c_j$ defined by the solution of  the linear equations system:

\begin{subequations}
\label{eq:discrete_variationalproblem}
\begin{align}
  \sum_{j=1}^{\ndof} \left(\int_{\Omega} \Delta \phi_j \Delta \phi_i \, \dx \right) c_j 
  &= \int_{\Omega} q \phi_i \, \dx, 
  && \text{for } i \in \mathcal{J}, \label{eq:dva} \\
  \sum_{j=1}^{\ndof} \phi_j(\xi_l) c_j 
  &= u_{\mathrm{e}_1}(\xi_l), \quad \xi_l \in \Gamma_{\mathrm{e}_1} , 
  && l = 1, \dotsc, n, \label{eq:dvc1} \\
  \sum_{j=1}^{\ndof} \partial_n \phi_j(\xi_l) c_j 
  &= u_{\mathrm{e}_2}(\xi_l), \quad \xi_l \in \Gamma_{\mathrm{e}_2} , 
  && l = 1, \dotsc, n \label{eq:dvc2}
\end{align}
\end{subequations}
or in more compact form
\begin{equation}\label{eq:Discrete_prob}
 \underbrace{\begin{pmatrix}
A_1 \\
C_1 \\
C_2
\end{pmatrix}}_{A}
c
=
\underbrace{\begin{pmatrix}
b^{\mathrm{I}}\\
b^{\mathrm{B}}
\end{pmatrix}\rule[-1.6em]{0pt}{1em}}_{b}
\end{equation}
with  
$A_1 \in \mathbb{R}^{|\mathcal{J}|\times\ndof}$,
$C_1 \in \mathbb{R}^{n \times \ndof}$, and
$C_2 \in \mathbb{R}^{n \times \ndof}$. 
The stiffness matrix $A_1$ is the discrete counterpart of the biharmonic differential operator.

We will first consider the quadratic case,  $n=(\ndof - |\mathcal{J}|)/2$, 
where the amount of data from inner points matches the missing amount 
of data on the boundary.   
By solving the above system (\ref{eq:Discrete_prob}) the required 
coefficients  $c = (c_j)_{j=1}^{\ndof}$ and hence $\uh$ are obtained. 
This will also  give the boundary data 
\[
f \vcentcolon = {\uh}_{\vert_{\partial \Omega _D}} \quad h:= { \partial _n \uh }_{\vert_{\partial \Omega _D}}
\]
resulting in the well-posed Dirichlet problem \cref{eq:ourproblem}. 

However, as already mentioned, the continuous problem 
\cref{eq:constrweekform} might have no solution and even in the case 
it has a unique solution the problem is ill-posed, \cite{lesnic1999numerical}. In the discrete case, this is reflected by an
ill-conditioned linear system, i.e., the matrix $A$ in 
\cref{eq:Discrete_prob} is near to a
rank-deficient matrix. Regularization techniques can be applied to make the solution more 
robust with respect to small data changes. They are based on 
introducing additional assumptions or  constraints.
Here we suggest the truncated singular value decomposition as a regularization method. Ambiguities
caused by poor or missing data are removed by asking for a minimal norm solution in a lower
dimensional subspace. Other regularization approaches can be found in 
\cite{hansen1994regularization} and the references therein.

We use the best rank-$k$ approximation $A_k$ of $A$ and compute the unknown coefficients by
\begin{equation} \label{eq:Truncated_svd_soln}
c_k  = A_k^{\dagger}b
\end{equation}
with $A_k^{\dagger}$ being the Moore--Penrose inverse of $A_k$. 

For this end we compute the singular value decomposition $A = U\Sigma V^T$ and replace its $\ndof-k$ smallest singular values
by zero to obtain $A_k$.   The regularization parameter $k$ must be 
determined. It controls the trade--off between the norm of the 
truncated solution and the norm of the corresponding residual:
\begin{equation}\label{eq:residual_vs_sol_norms}
|| c_k||_2^2  = \sum _{j=1}^k 
\frac{ (u_j\trans b)^2 }{\sigma _j^2}, \quad
||b - Ac_k||_2^2 =  \sum _{j=k+1}^{\ndof} (u_j\trans b)^2 
\end{equation}
Here,  $u_i$ is the $i$-th left singular vector corresponding to the singular value $\sigma_i$, 
see, e.g.~\cite{trefethen97}. 
An optimal truncation parameter $k$ can be obtained by using an L-curve 
method, which involves 
relating the two norms in \cref{eq:residual_vs_sol_norms} for varying 
values of $k$, \cite{aburidi2016comparative,hansen1999curve}.

\subsubsection{A synthetic example}
In the following example we illustrate the relation of the condition of 
the linear problem \cref{eq:Discrete_prob} to the distance of the 
inner curve $\Ge$ from the boundary $\partial \OD$. Additionally it demonstrates 
the relation between condition number and grid granularity.
  
\begin{figure}[h!]
\center{
  \includegraphics[width=0.5\textwidth]{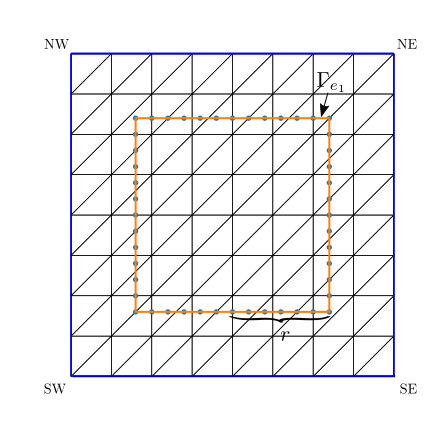}
}
\caption{Domain $\Omega$ with an inner curve $\Gamma_{e_1}$, forming a square of side length $2r$, discretized with 32 points. The background grid is $n_x \times n_x$ with $n_x = 8$. \label{Fig:squaredomain}}
\end{figure}

Consider the following problem

\begin{subequations} \label{synthProblem}
\begin{align}  
  \Delta ^2 u &= 0.1     &\text{in }         &\Omega = [0,2] \times 
  [0,2] \\
            u &= u_{e_1} &\text{on }  &\Gamma_{\mathrm{e}_1}  \\   
\partial _n u &= 0       &\text{on }  &\Gamma_{\mathrm{e}_2}  := \partial \Omega 
\end{align}
\end{subequations}
where $u_{e_1}$ is the solution of the corresponding  problem with 
homogeneous boundary conditions restricted to $ \Gamma_{\mathrm{e}_1}  $, the 
boundary of a centered inner
square with side length $2r < 2$, cf.~\cref{Fig:squaredomain}.

 This problem has a unique solution, but it is
ill-conditioned \cite{lesnic1999numerical}. The discrete solution space 
$\Vh$ is obtained by applying the virtual element method (VEM) using 
DUNE,
\cite{Dune2021}, on a triangular mesh constructed from an $n_x
\times n_x$ grid, $n_x$ being the number of intervals in each 
direction. Here, missing boundary data is
replaced by data from equally distributed points on the inner curve, see \cref{Fig:squaredomain}.
In the case of $n_x=8$ there are $\ndof=451$ degrees of freedom and the 
missing data is obtained from 32 points on the inner curve, \cref{Fig:squaredomain}.

\begin{figure}[hbt] 
  \includegraphics[width=0.55\textwidth]{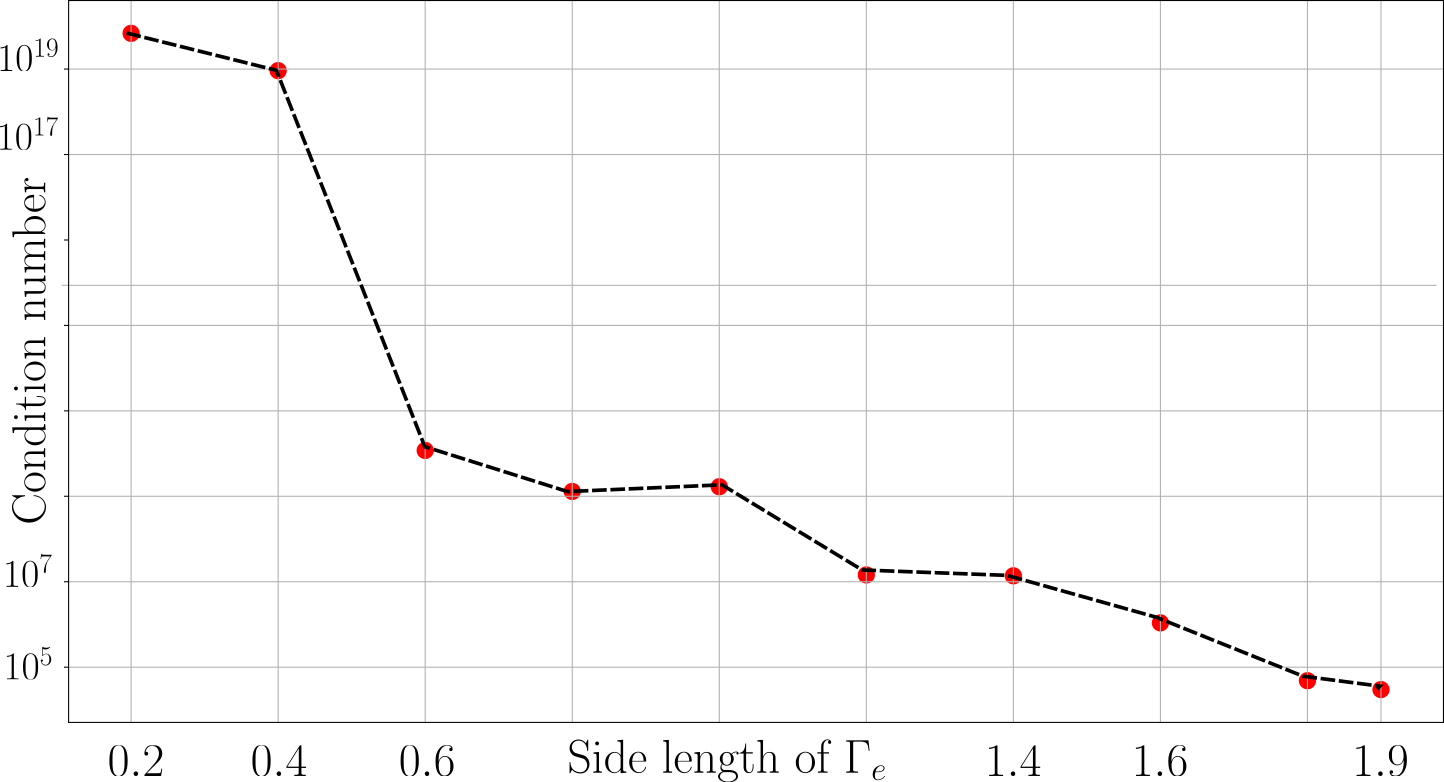} 
  \hfill
  \includegraphics[width=0.41\textwidth]{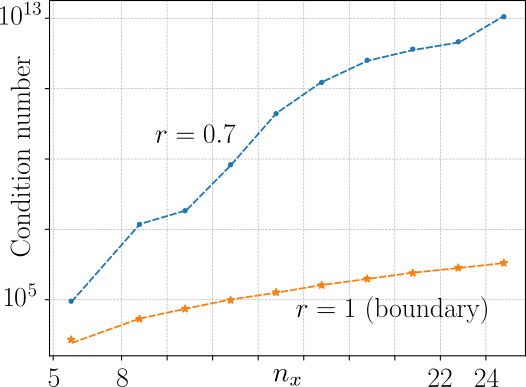}
\caption{{\it Left:} The condition number $\kappa(A)$ vs. the side length $2r$ of the inner
  curve. {\it Right:} The condition number $\kappa(A)$ vs. the number of intervals in each direction
  $n_x$ in the cases where $r=0.7$ and $r=1.0$.}
\label{Fig:syntheticProblem}
\end{figure}

Results are shown in \cref{Fig:syntheticProblem}. The picture to 
the left shows for $n_x=8$ how the
condition number depends on the distance between the outer boundary and the inner curve. As expected, the condition
number decreases as the side length of the inner curve increases and thus approaches the boundary.
The picture to the right shows how the condition number increases as the mesh is refined. The curves
show the condition numbers for different values of $n_x$ for the two cases: $r=0.7$ and
$r=1$, the latter is the solution of the biharmonic equation in the classical case, with the outer
boundary values known. The condition number in the case of $r=0.7$  
increases much faster with $n_x$. Especially, when the inner curve is 
far from the outer boundary, a regularization technique has to be 
applied when solving the linear system \cref{eq:Discrete_prob}.

\begin{figure}[h!]
\begin{center}
\includegraphics[width=0.4\textwidth]{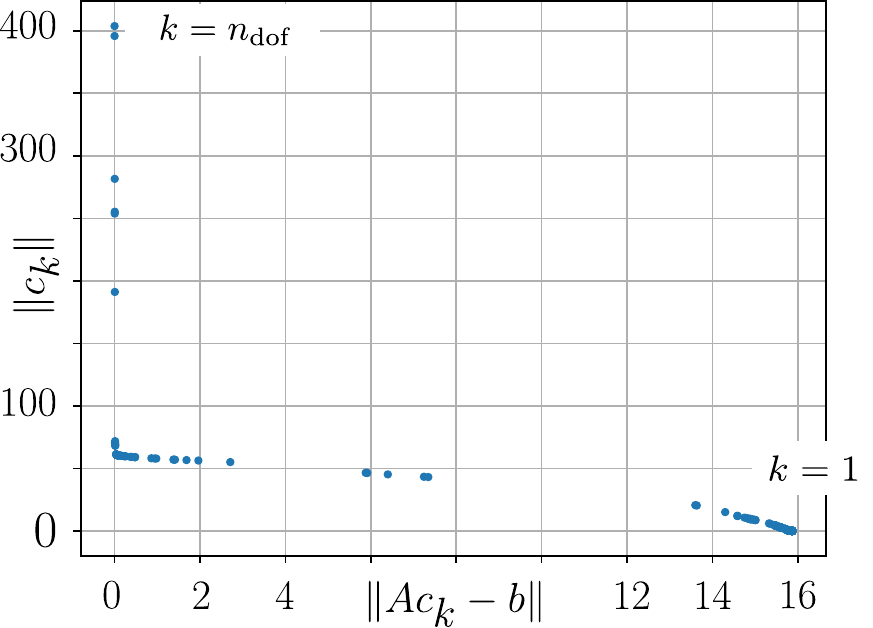}
\end{center}
\caption{Illustration of the norm of the truncated solution plotted 
against the norm of the resulting residual for truncation parameters  
$k \in \{1,2,\ldots,n_{\mathrm{dof}}=451\}$}.
\label{fig:L-curve}
\end{figure}

\Cref{fig:L-curve}  illustrates for problem \cref{synthProblem} the 
relation between the norms of solution and residual for different truncation 
parameters $k$ resulting in an L-shaped curve.

Here, one observes a reduction of the norm of the solution together 
with an increase of the norm of the residual when $k$ decreases.

\begin{figure}[hbt]
\center{
\includegraphics[width=0.5\textwidth]{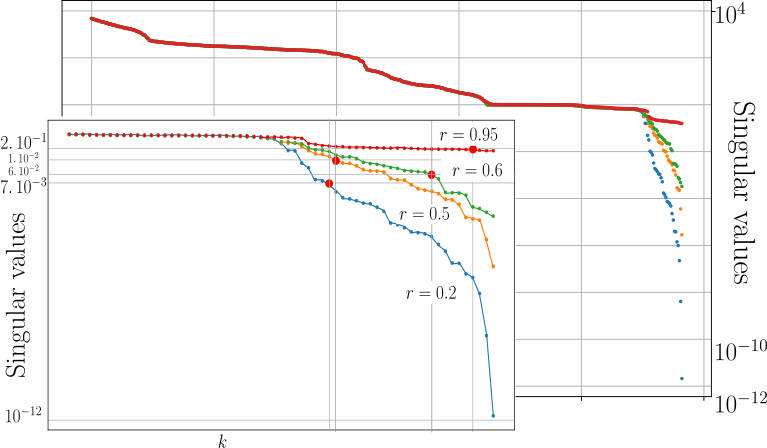}} \hfill
\includegraphics[width=.4\textwidth]{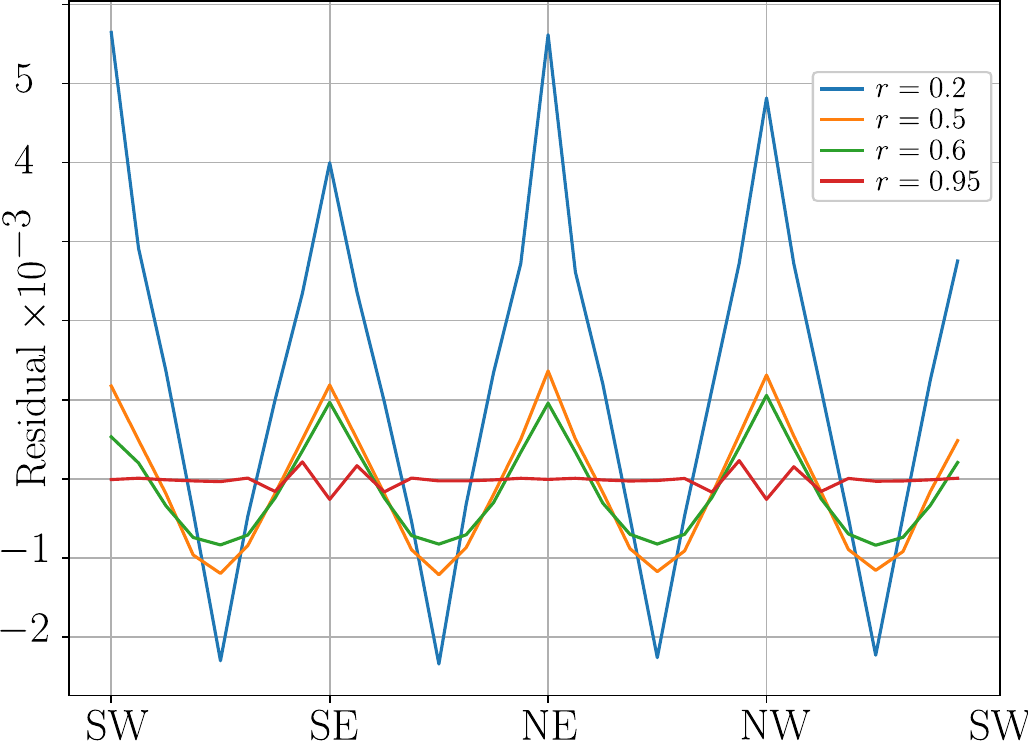}
\caption{{\it Left:} The singular values of the system matrix for square inner 
curves with different side lengths $2 r$. The embedded graph shows a 
zoom of relevant region with the smallest singular 
values and optimal truncation point in red. {\it Right:} Residual for 
different inner curves $\Gamma_{\mathrm{e}_1} $. It is the computed 
value on the boundary $\partial \Omega_{\mathrm{D}}$. \label{Fig:synteticSVD} 
}
\end{figure}

In \cref{Fig:synteticSVD} the singular 
values of $A$ are traced for different values of $r$. One can clearly observe that most singular values are largely unaffected by the distance of the inner curve from the boundary. However, beyond a certain value of $k$ they drop to different levels. The embedded figure focuses on this region  and 
indicates for each test case a choice of the regularization parameter 
$k$. The influence of the size of the inner square 
$\Gamma_{\mathrm{e}_1} $ on the computed the boundary data is shown in \cref{Fig:synteticSVD} 
(right). Recall that as $u_{e_1}$ was constructed from homogeneous 
boundary conditions, the figure shows the residual of the 
reconstruction experiment.  

\subsubsection{Case study: Modelling bedrock topography}
We return to the guiding example, presented in  \cref{sec:bedrock}, and 
compute the soil thickness by first computing the surface of the 
bedrock. The reason for this is that all subregions with exposed 
bedrock deliver data about the bedrock topography. That is the reason 
 why we here formulate the problem on the domain defined by the entire map region, $\Omega := \Om.$ For simplicity we assume that the second Dirichlet boundary 
condition is given by  $h = 0$, thus  $u_{\mathrm{e}_2} = 0$  and $ 
\Gamma_{\mathrm{e}_2} = \partial \Omega $ and $q =0$ in 
\cref{eq: var_pro}. Again, we aim for constructing the 
boundary function $f$. In the boundary segments intersecting with 
exposed bedrock, $\partial \Omega \cap \Oe$  this function is 
defined from the known bedrock topography, while on the remaining parts 
of the boundary we use available data from inner curves that is in this 
case topographic data from the boundaries of the 
exposed bedrock regions, $\partial \Oe$. 

The amount of unknowns on the boundary depends on the mesh. We 
can pick the same amount of data in a strategic way from these inner 
curves and an even more promising approach is to set up 
an overdetermined problem by picking more data than needed from all 
inner curves near the boundary, $\partial \Omega$. 

We demonstrate this here by considering an unstructured triangular mesh 
of the domain $\Omega$ generated by  Gmsh,~\cref{fig:meshes} 
(right). It consists of 
531 nodes, 1514 edges and 984 triangular elements. For setting up the finite-dimensional function space 
$\mathcal{V}^h$ we used an virtual element approach (VEM). The numbers 
of degrees of freedom is determined in this approach by the relation
\[
n_{\text{dof}} = n_{\text{nodes}} \cdot 3 + n_{\text{edges}} \cdot 
\left [ (\text{order} -3) + (\text{order} -2) \right ] + 
n_{\text{elements}} \cdot (\text{order} -3).
\]
We use a third order approach which results in $n_{\text{dof}} = 3107$ degrees 
of freedom. From these the desired 
surface is shaped. The values at all 76 boundary nodes are generated from inner curves 
$\Gamma_{e_1}$ that separate 
exposed from covered bedrock regions in the domain or any curve within the exposed region.  By taking $142$ 
points instead, the discrete problem \cref{eq:Discrete_prob} becomes overdetermined
with dimensions $n_{\mathrm{dof}}=3107$, $|\mathcal{J}|=2955$, $n=142$.

The singular values of $A$ range from $\sigma_{\max} = 20.64$ to 
$\sigma_{\min} = 5.14 \times 10^{-7}$. The L-curve method determines the 
truncation parameter $k= 3024$ that corresponds to the singular value 
$\sigma_k = 4.79\times 10^{-4}$ . 

\begin{figure}[hbt]
  \includegraphics[width=0.5\textwidth]{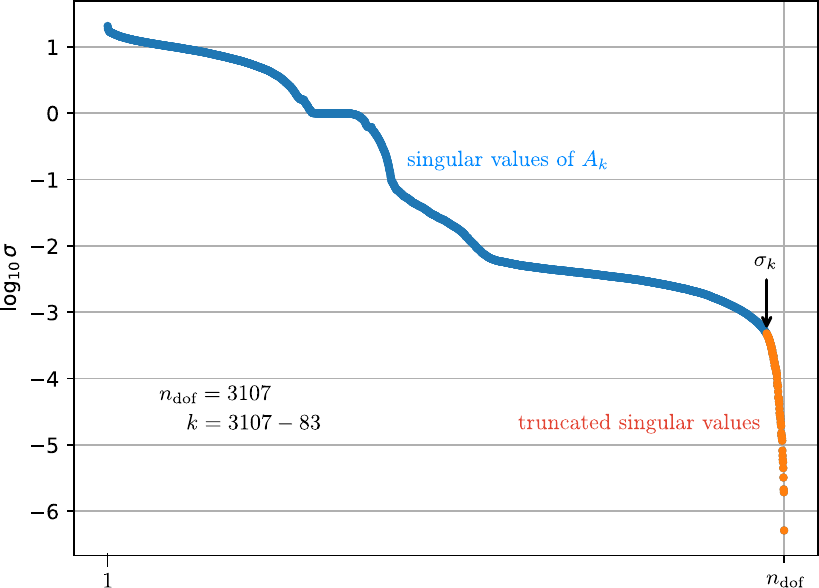}   
    \hfill     
  \includegraphics[width=0.4\textwidth]{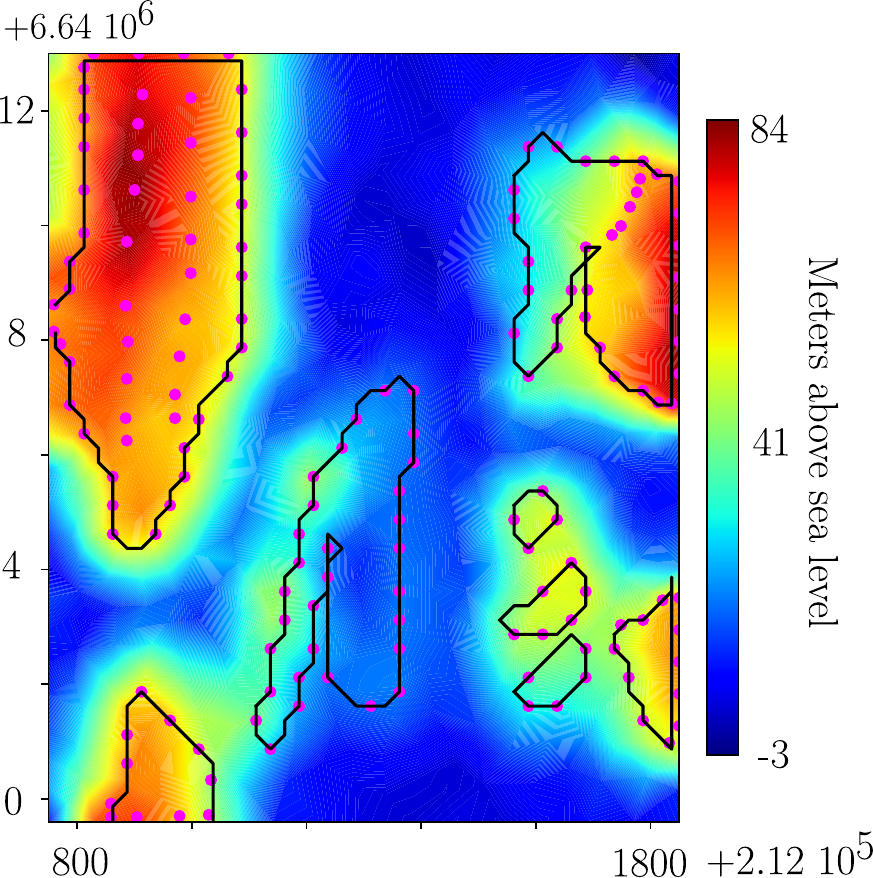}
    \caption{Bedrock topography discrete model: {\it Left:} Singular values of $A$ 
    and $A_k$. {\it Right:} Computed bedrock topography obtained by a 
    regularization based on truncating the 83 smallest singular values. The inner data points used are indicated by dots.\label{fig:L-shape}}
\end{figure}

By \cref{eq:Truncated_svd_soln} the minimal norm least 
squares solution is obtained, from which also the 
boundary values can be determined, \cref{fig:L-shape} (right).

Once the bedrock topography elevation is computed, we assess its physical realism by comparing it with the available terrain elevation data.

\begin{figure}[hbt]
  \begin{center}
  \includegraphics[width=.4\textwidth]{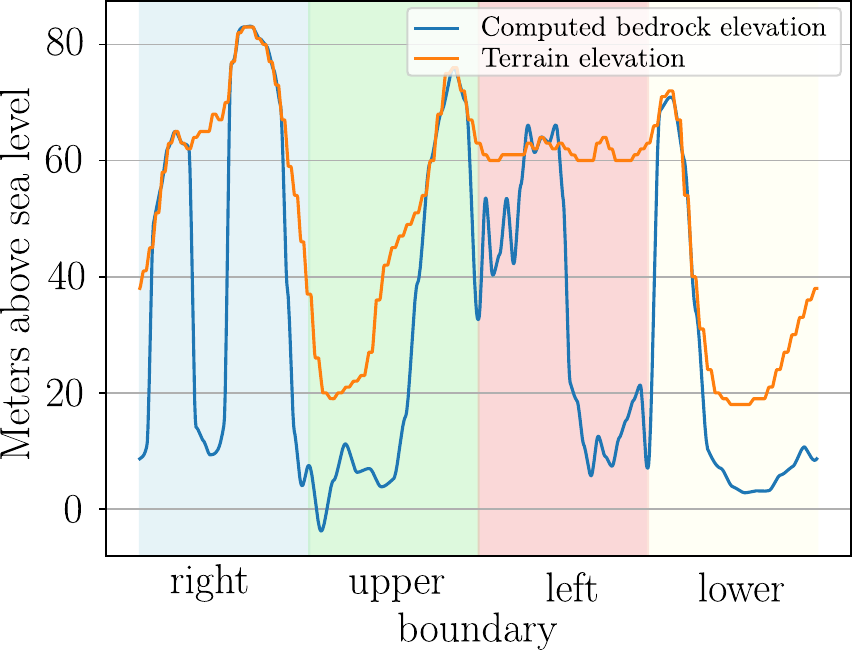} 
    \caption{Computed boundary values versus terrain elevations, cf. boundary in  \cref{fig:L-shape} (right) in counter clockwise direction. \label{Terrain vs Bedrock topo}}
  \end{center}
\end{figure}

As can be seen in Fig. \ref{Terrain vs Bedrock topo}, the bedrock topography on the boundary coincides with the terrain along exposed region portions of the boundary. In the rest of the places, as there are sediments, the bedrock topography lies strictly below the terrain topography.

\section{Software aspects} \label{sec:software}
In this section we describe the software tools used for computations in the soil thickness example.  
The approach of using the discrete biharmonic equation requires several practical steps. 
For most of them ready-to-use software is available.
In some steps a special adaptation of existing software is necessary.

The data acquisition step collects data describing the domain and measurements from public sources. 
The Python data analysis library \texttt{Pandas}, \cite{pandas} makes this data 
available for a domain description and meshing in Python. The domain is 
defined as a rectangular geographical region from which in \cref{sec:loadvector} the 
areas with exposed bedrock 
had to be excluded. The boundaries of these areas are obtained by 
determining the level set curves between the exposed region and the 
covered region by tools from the Python library \texttt{Matplotlib}. 
On this domain a triangular mesh is generated by \texttt{pygmsh}, a Python 
wrapper of \texttt{gmsh}, \cite{gmshsoftware}.

The mesh granularity highly affects computational complexity. It serves to 
define the surface modeling space $\Vh$. In contrast to the situation 
in PDE numerics the choice of the mesh granularity is not based on accuracy requirements 
when seeking a solution of a PDE. It is rather a design parameter for 
shaping the surface. This is essential when deciding on numerical 
algorithms as especially the algorithms from linear algebra used in 
this context have 
a computational complexity that grows cubically with the dimension.

To formulate the biharmonic equation on that domain and to set up the 
discrete formulation of \cref{eq:ourproblem} the modular toolbox for solving 
partial differential equations \texttt{DUNE} \cite{Dune2021} is used. 
It has tools to interface with the mesh data and a mesh-to-DGF 
converter to mark the different boundary segments.

Our intention is to generate a smooth surface by $\mathcal{C}^1$ 
functions. For this end the DUNE module \texttt{dune-vem} is used to
set up the discrete function spaces by a third order virtual element approach,  \cite{brezzi2013virtual,beirao2013basic}. 
This method is an off-spring of the classical finite element method, 
that allows to solve the biharmonic equation but with less degrees 
of freedom than a finite element approach with the same order would have 
needed.
Thus, the VEM generated matrices are of smaller dimension, which leads to a higher computational efficiency.

In \cref{sec:boundary}, DUNE is used to generate the system matrices and 
vectors in \cref{eq:Discrete_prob}. From that, the solution of the 
discrete problem is computed outside of DUNE by the use of the Python module 
\texttt{scipy.linalg}. 
 
The results are visualized by graphical tools 
provided by \texttt{DUNE} and the Python module \texttt{matplotlib}.


\section{Summary}
In this paper strategies for construction surfaces as solutions of the discrete biharmonic equation
based on a rather small amount of data have been presented. These strategies result in methods
which have been implemented based on standard free software. The ideas are quite general, and
the solutions depend on choice of load functions, parameters and the choice of data points on inner
curves, which we believe make them applicable to a larger set of situations. In this paper the
strategies have 
been demonstrated on a concrete problem from soil science to construct surfaces mimicking some trends in
soil thickness locally on covered areas, as well as the bedrock topography over a larger domain.
The results are very useful as basis for further refinement of the physical results by use of
geostatistical means.


\section*{Acknowledgement}
We want to acknowledge the opportunity to have many inspiring discussions 
with our colleagues Markus Grasmair, Robert Klöfkorn, Andreas Dedner, 
Alemayehu Adugna, and
Zerihun Knife.

%
\section*{Conflict of interest}
The authors declare that they have no conflict of interest.

\bibliographystyle{spmpsci}      
\bibliography{bibfile.bib}   

\end{document}